\documentclass[12pt]{article}
\usepackage{amsfonts,amssymb,amsmath,bm}
\usepackage{a4,graphics,epsfig,psfrag,geometry}
\usepackage[utf8]{inputenc}
\usepackage{psfrag}
\usepackage{xcolor}
\usepackage{tabularx}
\usepackage{float}
\newcolumntype{Y}{>{\raggedright\arraybackslash}X}
\newcommand{\RR}{\mathbb R}

\newcommand{\f}[1]{\mathbf{#1}}

\begin{document}
\begin{center}
{\LARGE \bf 
Isogeometric Simulation of Thermal Expansion
for Twin Screw Compressors} \\[4mm]
A Shamanskiy$^1$ and B Simeon$^1$ \\[1mm]
$^1$ Dept. of Mathematics, TU Kaiserslautern \\
Paul-Ehrlich-Str.~31, Kaiserslautern, Germany \\ 
E-mail: \{shamansk$|$simeon\}@mathematik.uni-kl.de  
\end{center}

\begin{quote} \small
	
{\bf Abstract:} 
Isogeometric Analysis (IGA) is a recently introduced computational approach intended to breach the gap between the Finite Element Analysis and the Computer Aided Design worlds. In this work, we apply it to numerically simulate thermal expansion of oil free twin screw compressors in operation. High global smoothness of IGA leads to a more accurate representation of the compressor geometry. We utilize standard tri-variate B-splines to parametrize the rotors, while the casing is modeled exactly by using NURBS. We employ the Galerkin version of IGA to solve the thermal expansion problem in the stationary case. The results allow to estimate the contraction of the clearance space between the casing and the rotors. The implementation is based on the open source C++ library G+Smo.

This work is supported by the European Union within the Project MOTOR:
Multi-ObjecTive design Optimization of fluid eneRgy machines.

Thermo-mechanical coupling.
\end{quote}

\section{Introduction}
A twin screw compressor is a type of gas compressor with a great practical significance in various branches of engineering~\cite{stosic2005screw}. In this work, we focus on the oil free subclass of these machines. The absence of oil and its cooling effect leads to a much higher temperature rise than in oil injected compressors. The resulting thermal stress distorts the shape of the rotors, making it necessary to increase clearances in order to avoid contact between the rotors. By considering the effects of thermal expansion in numerical simulations~\cite{spille2015cfd}, an estimate of the minimal safe distance can be acquired and the efficiency of compressors could be increased.

Here we propose an alternative approach to thermal expansion simulation with the use of Isogeometric Analysis (IGA), which is a modern computational method breaching the gap between the Finite Element Analysis and the Computer Aided Design worlds~\cite{Cottrell.2009,Hughes2005}. NURBS-based IGA is particularly attractive for applications with rotor-type geometries and where spinning motion is involved, due to its ability to exactly represent circular arcs. Another key feature of IGA is its high global geometric smoothness which allows to achieve greater accuracy per degree of freedom when modeling intricate shapes. In this paper, we demonstrate one possible strategy of generating 3D IGA models of the rotors and the casing. We start from the point cloud representation of the rotors' profiles and approximate it with cubic B-spline curves. Using these approximations, the 3D tensor product B-spline parametrizations of the rotors are generated, while the casing model is built with NURBS. We then state a mathematical model of thermal expansion and use it to simulate thermal growth in the case of stationary rotors. The results are in a good agreement with the previous work, but much less computational resources are needed. We then use the results to estimate the contraction of the clearance space between the casing and the rotors.

We give now a short overview of different simulation methods and aspects for twin screw compressors. One of the most common approaches, 0D chamber-model simulations, is described in the work by Kauder~\cite{kauder2003adiabatic}. For a systematic introduction to 3D CFD simulations we refer to Kovacevic et al.~\cite{kovacevic2007screw}, while advanced CFD simulations including leakage flows are presented in ~\cite{pascu2014numerical,andres2016cfd}. A novel efficient tool for fluid domain meshing is proposed by Hesse~\cite{hesse2014structured}. In a recent work by Utri~\cite{utri2017energy} the optimization of screw machines with variable rotor pitch is addressed.

This paper is organized as follows. In Section 2, we describe the creation of geometric models for the rotors and the casing. Section 3 deals with the mathematical model of thermal expansion and demonstrates its use on a numerical example. In Section 4, we draw conclusions and outline current and future research directions.

\section{Geometric Models}
\begin{figure}[H]
	\centering
	\includegraphics[height = 5.5cm]{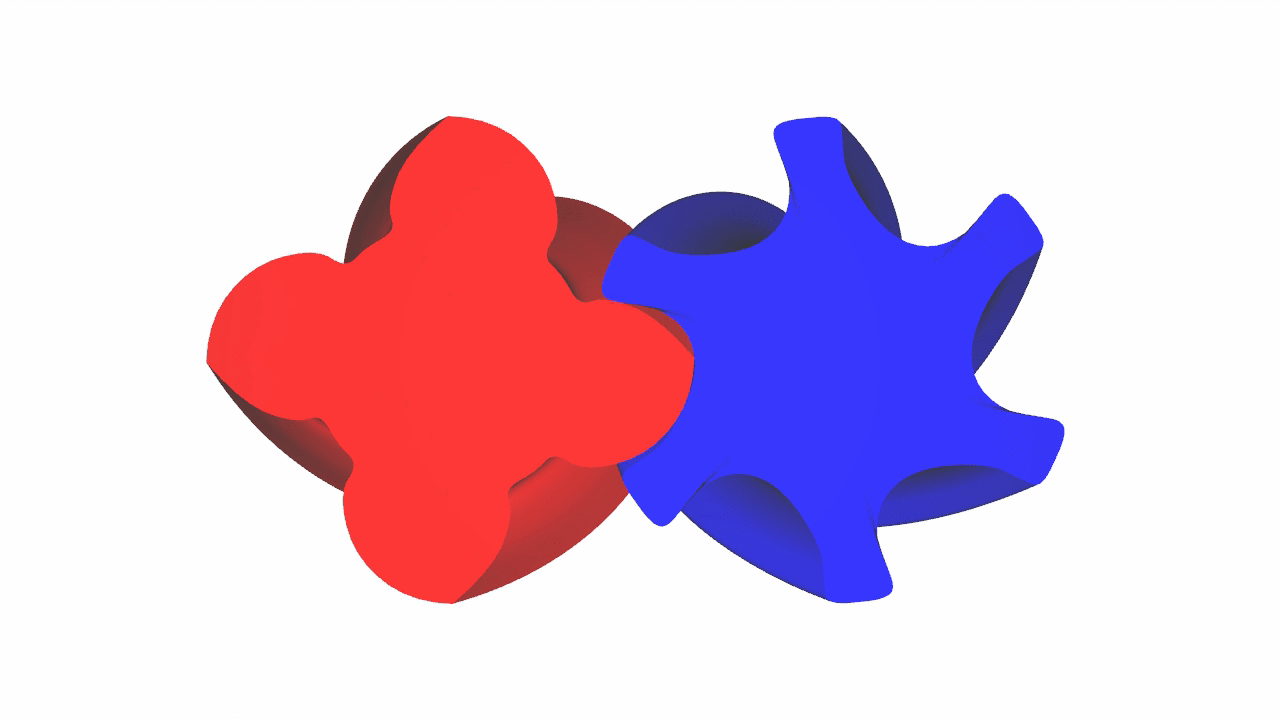}
	\caption{Screw machine's rotors. Throghout this article they are refered to as the male (left, red) and the female (right, blue).}\label{fig:rotors}
\end{figure}

Both of the compressor rotors possess an intricate geometry (Fig.~\ref{fig:rotors}), a careful representation of which is the focal point of this article. An extremely small clearance height between the rotors and the surrounding casing requires a high degree of accuracy from the approach chosen to describe these objects. The conventional approach is the use of point clouds. To run numerical simulations, however, another representation is necessary. Within the framework of Isogeometric Analysis (IGA) the well-known tools from the world of Computer-Aided Design (CAD) - B-splines or NURBS - are utilized to model objects. However, unlike in CAD, not only the boundaries or surfaces of objects are modeled, but also their interior. The standard way to do this in IGA are tensor product B-splines.

\subsection{Tensor Product B-splines}
To fix the notation, we shortly outline the construction of tensor product B-splines
in the case of $d=2$ dimensions. First we specify the 
polynomial degrees $(p,q)$ and the horizontal and vertical knot vectors 
\[
\Xi = \left\{\xi_1\le \xi_2\le\ldots \le \xi_{m+p+1}\right\},\quad 
\Psi = \left\{\eta_1\le \eta_2\le\ldots \le \eta_{n+q+1}\right\},
\]
which contain non-decreasing parametric real values so that
\[
0\le\mu(\Xi,\xi)\le p+1 \qquad \text{and} \quad
0\le\mu(\Psi,\eta)\le q+1 
\]
are the multiplicities of the parameter values in the knot vectors 
(the multiplicity $\mu(X,x)$ is zero if the given value $x$ is not a knot in $X$).
We write $M_{i,p}(\xi)$  and $N_{j,q}(\eta)$ for the univariate B-splines, 
which are computed by means of the recursion \cite{Boor1978}
\begin{align} 
N_{j,0}(\eta) &= \left\{ \begin{array}{ll}
1 & \text{for } \eta_j\le \eta < \eta_{j+1},\\
0 & \text{otherwise},
\end{array} \right.
\label{eq:bspline1}\\
N_{j,q}(\eta) &= 
\frac{\eta-\eta_j}{\eta_{j+q}-\eta_j}N_{j,q-1}(\eta) +
\frac{\eta_{j+q+1}-\eta}{\eta_{j+q+1}-\eta_{j+1}}N_{j+1,q-1}(\eta)
\label{eq:bspline2}
\end{align}
for $j=1,\ldots, n$. Fractions with zero denominators are considered zero. The same recursion applies to $M_{i,p}(\xi)$, replacing $j,q$ in 
(\ref{eq:bspline1})--(\ref{eq:bspline2}) by $i,p$ for $i=1,\ldots, m$. 

A function $\f{f}(\xi,\eta): 
[\xi_{p+1},\xi_{m+1}]\times[\eta_{q+1},\eta_{n+1}] \rightarrow \mathbb{R}^2$ is called a bivariate tensor product B-spline function if 
it has the form
\begin{equation}\label{eq:bsplinef}
\f{f}(\xi,\eta)=\sum_{i=1}^{m}\sum_{j=1}^{n}  M_{i,p}(\xi)N_{j,q}(\eta)  \f{d}_{i,j} 
\end{equation}
with the de Boor control points $\f{d}_{i,j} \in \mathbb{R}^2$, which form the control net
associated to the parametric representation. In order to work with the unit square as 
parameter domain, combined with open knot vectors, we assume
additionally $\xi_1=\ldots =\xi_{p+1} = 0$ and 
$\xi_{m+1} =\ldots = \xi_{m+p+1} = 1$ 
at the beginning and at the end; analogously  $\eta_1=\ldots =\eta_{q+1} = 0$ and 
$\eta_{n+1} = \eta_{n+2}=\ldots = \eta_{n+q+1} = 1$.

In $d=3$ dimensions, a trivariate B-spline tensor product function is generated
analogously, which results in 
\begin{equation}\label{eq:bsplinef3d}
\f{f}(\xi,\eta,\psi)=\sum_{i=1}^{m}\sum_{j=1}^{n}\sum_{k=1}^\ell  M_{i,p}(\xi)N_{j,q}(\eta)L_{k,r}(\psi)  \f{d}_{i,j,k} 
\end{equation}
with an additional set of univariate B-splines $L_{k,r}$ and control 
points $\f{d}_{i,j,k} \in \mathbb{R}^3$.

\subsection{Model Construction}
\begin{figure}[H]
	\centering
	\includegraphics[height = 6cm]{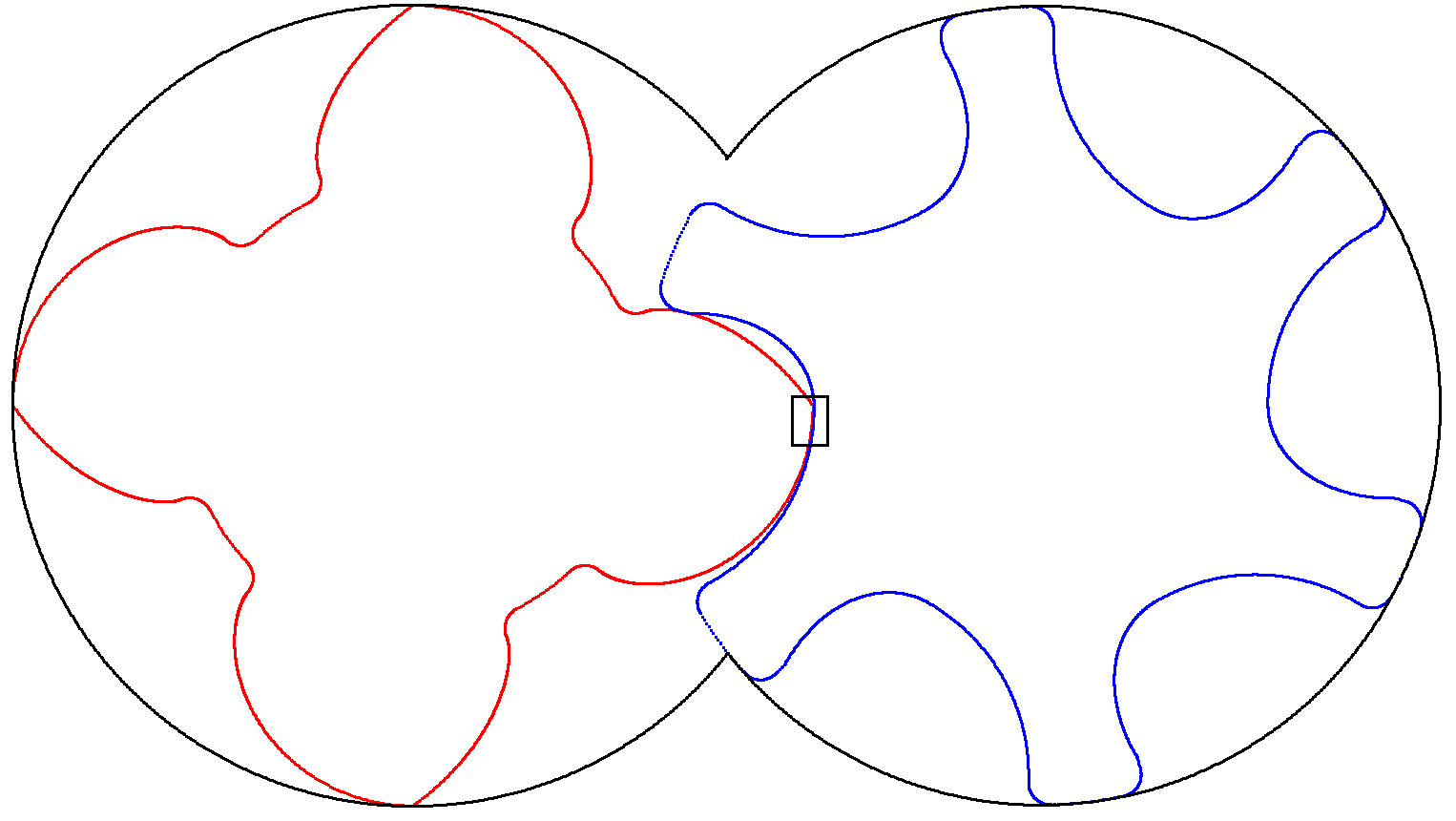}
	\includegraphics[height = 6cm]{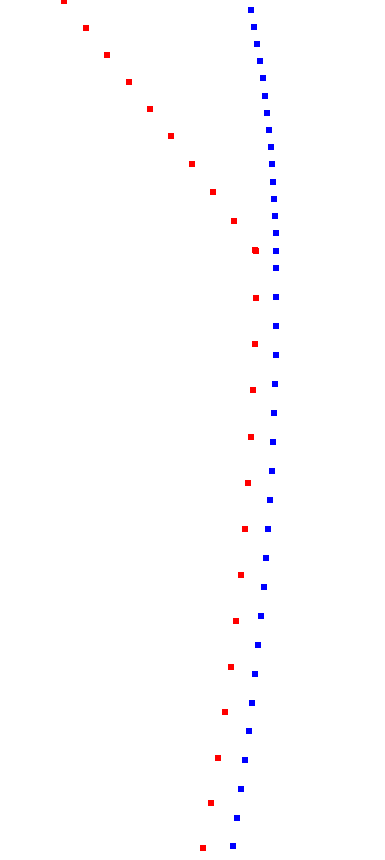}
	\caption{Point cloud describing the male and the female rotors as well as the surrounding casing. Magnified segment to the right clearly depicts the points and the clearance height.}\label{fig:cloud}
\end{figure}

The geometry was originally provided in the form of a point cloud where the male and the female rotors of SRM type are described with 2572 and 2292 points correspondingly, see Fig.~\ref{fig:cloud}. It was scaled to match the configuration from~\cite{spille2015cfd}. As a result, the axis distance is 80 mm, the clearance height is 44 $\mu$m between the casing and the rotors  and approximately 100 $\mu$m between the rotors, the radius of the rotors is 50.97 mm. The difference in scales makes it challenging to come up with a B-spline approximation of rotors' profiles without a large number of control points. That is one of the reasons why point clouds, rather than CAD models, are usually used to describe these objects.

We consider a B-spline approximation accurate enough if it deviates from the point cloud for less than 10\% of the clearance height. This accuracy can be achieved by using cubic splines with 73 control points for the male rotor and with 109 control points for the female rotors. Four sharp corners of the male rotor are modeled by increasing the knots multiplicity. For the technical details of B-spline curve approximation we refer to \cite{Piegl1997}. Two circular arcs forming the casing can be modeled exactly by NURBS. This ability to exactly represent circular shapes makes IGA particularly suitable for applications where spinning motion is involved.

\begin{figure}
	\centering
	\includegraphics[height = 6.9cm]{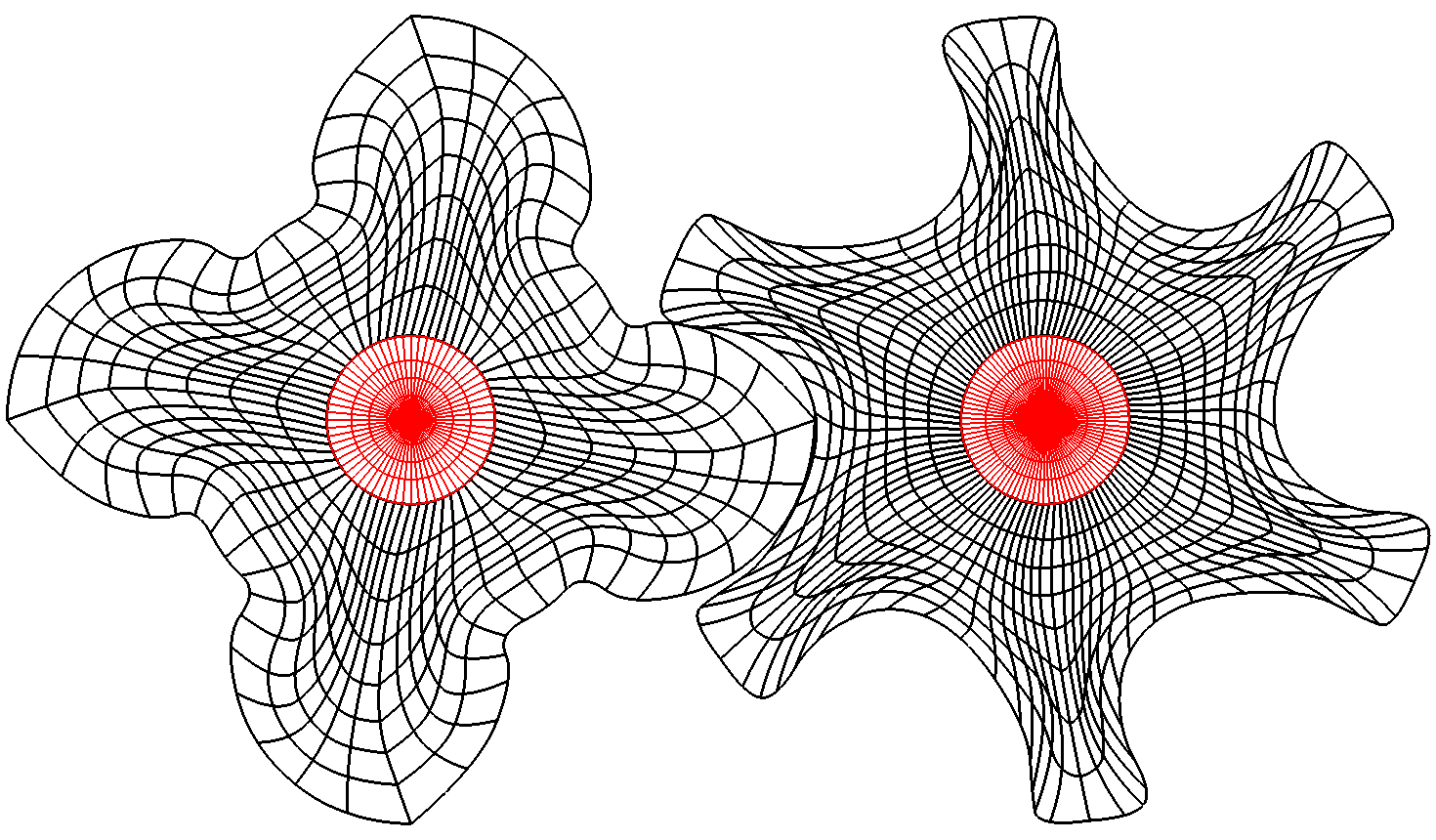}
	\includegraphics[height = 6.9cm]{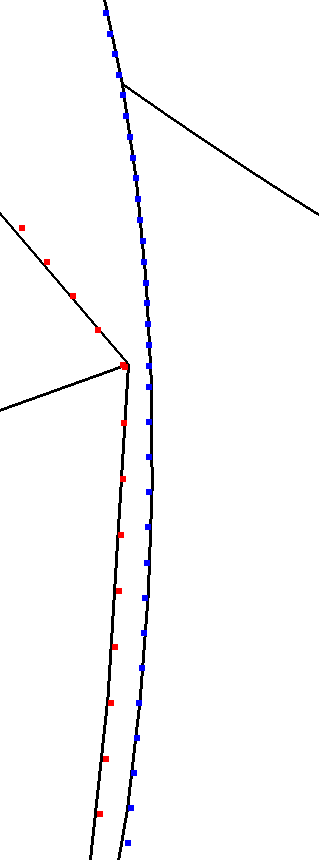}
	\caption{2D scaled boundary parametrizations for the rotors. Original point cloud and its B-spline approximation (to the right, magnified).}\label{fig:meshes2d}
\end{figure}

Next we generate 2D parametrizations for the cross-sections of the rotors. We choose the Scaled Boundary approach~\cite{Arioli2017} which is natural for rotary objects and allows to describe the interior of each rotor with a minimum number of patches, see Fig.~\ref{fig:meshes2d}. The outer layer is formed by "scaling" the boundary towards the circular patch in the center. This central patch is later used to construct a shaft. The parametrization is $C^2$-continuous inside the patches, but only $C^0$ between them. Both patches are obtained as a geometric mapping from a unit square. For the central patch this mapping is not bijective at the center. However, despite this singularity, this type of parametrizations is analysis-suitable. 

\begin{figure}[H]
	\centering
	\includegraphics[height = 6.9cm]{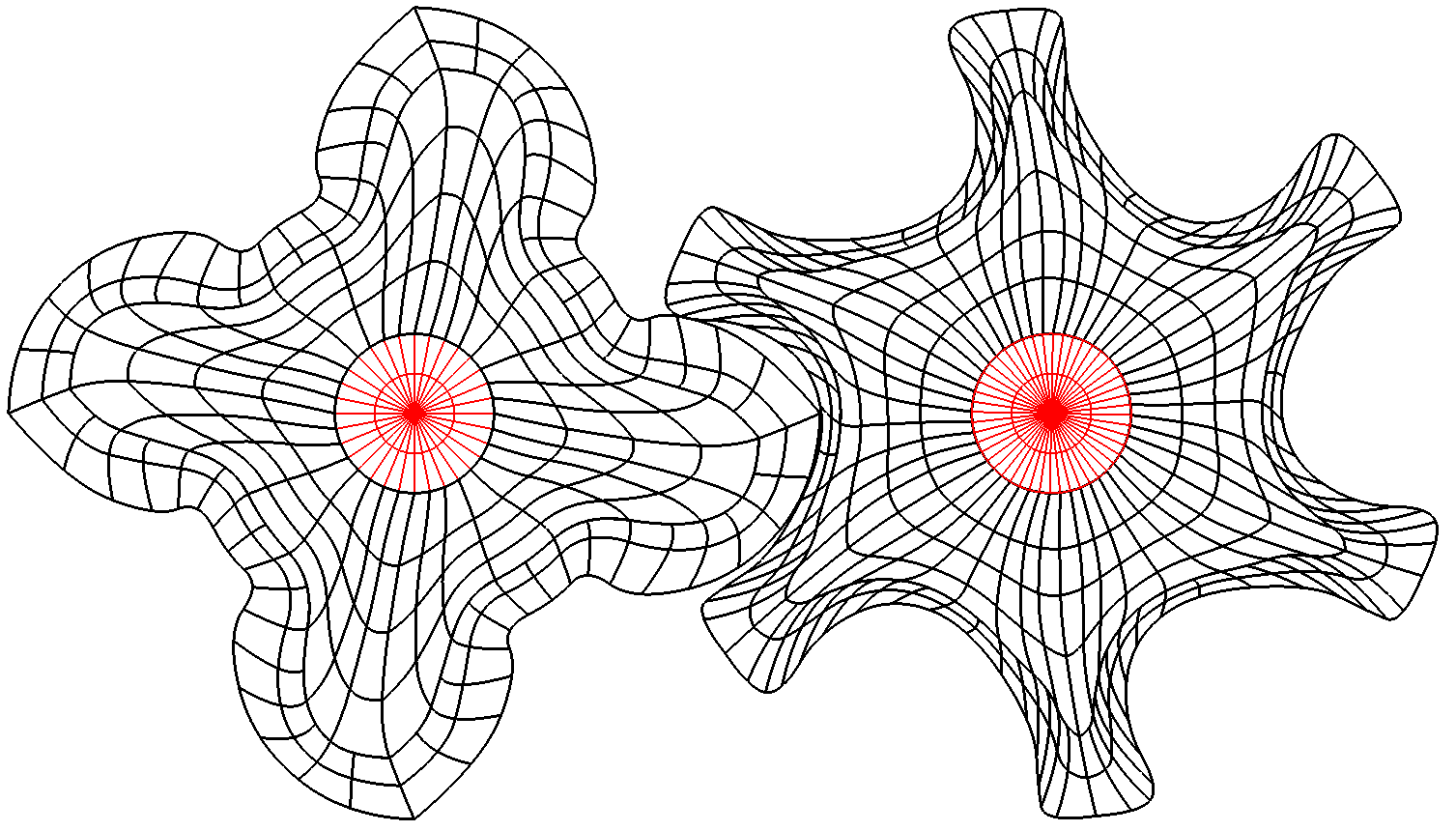}
	\caption{Local refinement with THB-splines allows to allocate degrees of freedom denser along the boundaries of the rotors than in their interior.}\label{fig:thb}
\end{figure}

Optionally, the planar parametrizations can be locally refined or coarsened using THB-splines~\cite{vuong2011hierarchical}. By locally refining boundary layers of the rotors a finer resolution can be achieved there. This is especially attractive for the twin screw compressor example since the surfaces of the rotors require a finer resolution for analyzing and optimizing the mechanical properties. Figure~\ref{fig:thb} depicts a locally coarsened planar mesh. The B-spline approximation of the boundaries is maintained while less degrees of freedom are used in the interior of the rotors.

\begin{figure}
	\centering
	\includegraphics[height = 6cm]{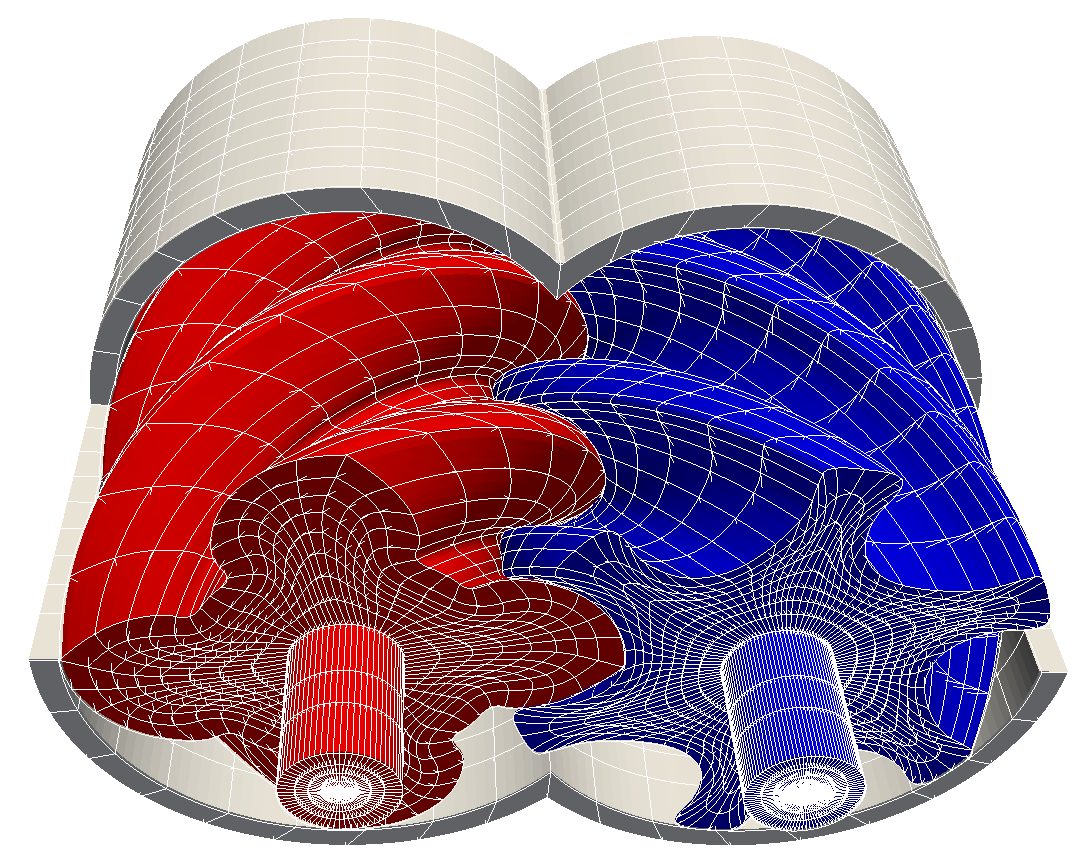}
	\includegraphics[height = 4cm]{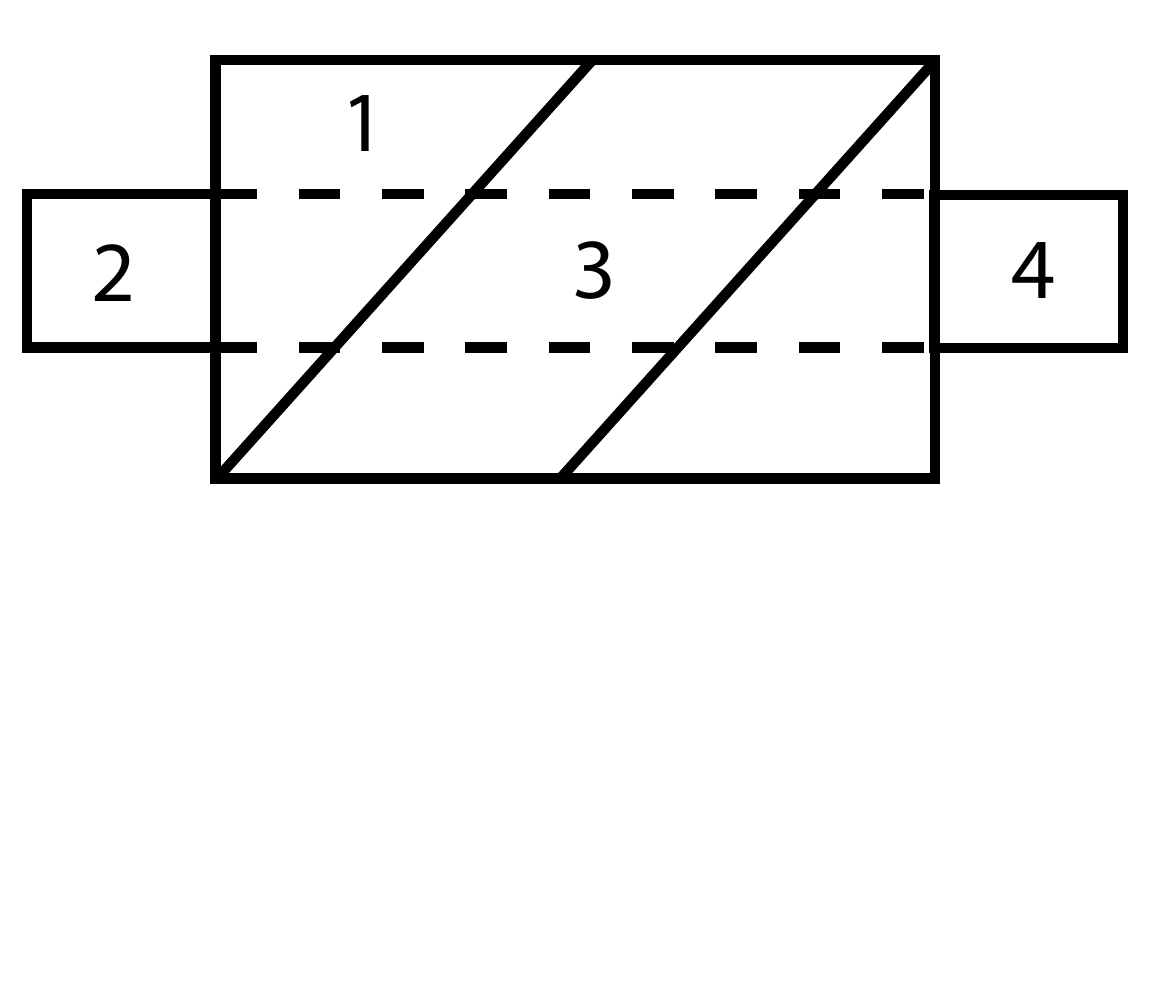}
	\caption{3D models of the rotors and the casing. Multi-patch structure of the parametrizations: each rotor consists of 4 patches with matching boundary descriptions.}\label{fig:rotor3d}
\end{figure}

The final step is a 3D model. Each rotor is represented by a multi-patch geometry with 4 patches, see Fig.~\ref{fig:rotor3d}. The shaft has to be split into three parts so that the central part could have the matching boundary description with the inner surface of the outer layer. The control points for the outer layer are obtained by combining several layers of rotated control points from the 2D cross-section. 20 layers of control points are used for both rotors. Both rotors have a length of 168.3 mm, the male is modeled with a pitch of 300$^\circ$, the female with -200$^\circ$. The shaft patches 2 and 4 are 30 mm long and the casing is 6 mm thick.

\section{Thermal Expansion}
We consider here a very basic model of thermal expansion, which is nevertheless applicable for the case of screw compressors. It is based on the equations of linear elasticity, where the thermal growth is represented by an extra term. The use of a linear model is justified since only small deformations appear during the operation of the compressor, otherwise the machine would be blocked. However, these small deformations may still alter the shape of the clearing and change the efficiency of the compressor.

\subsection{Mathematical Model}
A body of isotropic material subject to a heat change from an initial temperature $T_0$ to a temperature $T$ expands or contracts uniformly in all directions~\cite{Gross2011}. This can be modeled by the thermal strain
\begin{equation}
\label{eq:thstr}
\pmb{\varepsilon}_T = \alpha (T-T_0)\f{I},
\end{equation}
where $\alpha$ is the coefficient of thermal expansion of the material and $\f{I}$ is the $d\times d$ identity matrix. The resulting deformation is described by the displacement field $\f{u}\in \mathbb{R}^d$. To take into account the mechanical load as well we consider the linearized mechanical strain $\pmb{\varepsilon}(\f{u}):=\nabla^{(s)}\f{u}$, where $\nabla^{(s)}$ is the symmetric gradient. We define the total strain $\pmb{\varepsilon}_S(\f{u})$ as a superposition of both:

\begin{equation}
\label{eq:super}
\pmb{\varepsilon}_S(\f{u}) = \pmb{\varepsilon}(\f{u}) + \pmb{\varepsilon}_T = \frac{1}{2}(\nabla\f{u}+\nabla\f{u}^T) + \pmb{\varepsilon}_T.
\end{equation}

The total stress $\pmb{\sigma}(\pmb{\varepsilon}_S(\f{u}))$ is related to the total strain $\pmb{\varepsilon}_S(\f{u})$ by Hooke's law with the Lam\'{e} parameters $\lambda$ and $\mu$: 
\begin{equation}
\label{eq:hooke}
\pmb{\sigma}(\pmb{\varepsilon}_S(\f{u})) = \f{C}\pmb{\varepsilon}_S(\f{u}) = 2\mu\pmb{\varepsilon}_S(\f{u}) + \lambda\text{tr}(\pmb{\varepsilon}_S(\f{u}))\f{I}.
\end{equation}

With the total stress defined we can now write the system of linear elasticity equations for the domain $\Omega\subset\mathbb{R}^d$ as

\begin{align}\label{eq:equil1}
-\operatorname{div} \pmb{\sigma}(\pmb{\varepsilon}_S(\f{u})) &= \f{g} \text{ in } \Omega,\\
\f{u} &= \f{0} \text{ on } \partial\Omega_D,\\
\label{eq:equil2}
\pmb{\sigma}(\pmb{\varepsilon}_S(\f{u}))\cdot \f{n} &= \f{f} \text{ on } \partial\Omega_N.
\end{align}

Since both the divergence operator and Hooke's law are linear, we can move the thermal terms to the right-hand side of the equations (\ref{eq:equil1}) and (\ref{eq:equil2}). Keeping in mind a change of the sign for the outer normal, we rewrite the equations in the following form:

\begin{align}
\label{eq:elthermo1}
-\operatorname{div}\pmb{\sigma}(\pmb{\varepsilon}(\f{u})) &= \f{g} - \operatorname{div}\pmb{\sigma}(\pmb{\varepsilon}_T) \text{ in } \Omega,\nonumber\\
\f{u} &= \f{0} \text{ on } \partial\Omega_D,\\
\pmb{\sigma}(\pmb{\varepsilon}(\f{u}))\cdot \f{n} &= \f{f} + \pmb{\sigma}(\pmb{\varepsilon}_T)\cdot \f{n} \text{ on } \partial\Omega_N.\nonumber
\end{align}

Now for any given temperature distribution $T$ the displacement field $\f{u}$ can be computed. The temperature itself is assumed to be given or can be acquired as a solution of Poisson's equation or the heat equation depending on whether the stationary or the transient problem is considered. We assume that the temperature change is small enough so that the material properties remain constant. Another assumption is that the domain deformation does not affect the temperature distribution inside. Then the problem is only one-way coupled with the temperature change driving the thermal expansion without feedback. 

\subsection{Galerkin IGA}
Let $\mathcal{V} = \{\f{f}\in H^1(\Omega)^d \text{ }|\text{ } \f{f} = \f{0}\text{ on } \partial\Omega_D\}$ be the trial solution and the weighting function spaces, which are assumed to be the same. Then the weak form of the system (\ref{eq:elthermo1}) is defined in the following way, which is equivalent to the Principle of Virtual Work:
\begin{align}
\label{eq:weakform}
Find\,\f{u} \in \mathcal{V},\,such\,that&\nonumber\\
\int_\Omega\nabla^{(s)}\f{u}:\f{C}\,\nabla^{(s)}\f{v}\,\text{d}\f{x} &= \int_\Omega(\f{g}-\operatorname{div}\pmb{\sigma}(\pmb{\varepsilon}_T))\cdot\f{v}\,\text{d}\f{x}\\
&+ \int_{\partial\Omega_N}(\f{f}+\pmb{\sigma}(\pmb{\varepsilon}_T)\cdot\f{n})\cdot\f{v}\,\text{d}\f{s}, \quad for\,\forall \,\f{v} \in \mathcal{V}.\nonumber
\end{align}
Here : denotes the Frobenius inner product of two matrices. For functions $\f{u},\f{v}\in\mathcal{V}$, the bilinear form 
\begin{equation}
a(\f{u},\f{v}) = \int_\Omega\nabla^{(s)}\f{u}:\f{C}\,\nabla^{(s)}\f{v}\,\text{d}\f{x}
\end{equation}
is well-defined, and even more, it is symmetric and coercive. Setting 
\begin{equation}
l(\f{v})= \int_\Omega(\f{g}-\operatorname{div}\pmb{\sigma}(\pmb{\varepsilon}_T))\cdot\f{v}\,\text{d}\f{x} + \int_{\partial\Omega_N}(\f{f}+\pmb{\sigma}(\pmb{\varepsilon}_T)\cdot\f{n})\cdot\f{v}\,\text{d}\f{s}
\end{equation}
as linear form for the integration of the right hand side, the equation of the weak form (\ref{eq:weakform}) can be written as
\begin{equation}\label{eq:weakform2}
a(\f{u},\f{v}) = l(\f{v}).
\end{equation}
The Galerkin projection of the weak form (\ref{eq:weakform}) 
in physical coordinates $\f{x}$
replaces the infinite dimensional space $\mathcal{V}$ by a finite dimensional subspace $\mathcal{V}_h \subset \mathcal{V}$, with the subscript $h$ indicating the relation to a spatial grid. 
Let $\pmb{\phi}_1, \ldots, \pmb{\phi}_K$ be a basis of $\mathcal{V}_h$, then the numerical approximation $\f{u}_h$ is constructed as linear combination
\begin{equation}\label{num:defuh}
\f{u}_h =  \sum_{i=1}^{K} q_i \pmb{\phi}_i
\end{equation}
with unknown real coefficients $q_i \in \RR$. 

Upon inserting $\f{u}_h$ into the weak form (\ref{eq:weakform2}) and testing with $\f{v} = \pmb{\phi}_j$ for $j=1,\ldots,K$, one obtains the linear system
\begin{equation}\label{num:linsys}
\f{A} \f{q} = \f{r}
\end{equation}
with $K \times K$ stiffness matrix $\f{A} = ( a(\pmb{\phi}_i,\pmb{\phi}_j))_{i,j=1,\ldots,K}$ and right hand side vector $\f{r} = (l( \pmb{\phi}_i))_{i=1,\ldots,K}$.
Since the matrix $\f{A}$ inherits the properties of the bilinear form $a$, it is straightforward to show that $\f{A}$ is symmetric positive definite, and thus the numerical solution $\f{q}$ or $\f{u}_h$, respectively, is well-defined.

Let $\mathcal{S}$ denote the function space of tensor product B-splines defined in Section 2. In Galerkin-based IGA, we define 
$\pmb{\phi}_i = \pmb{\psi}_i \circ \f{F}^{-1}, \, i=1,\ldots,K,$  
as basis functions
via the push forward operator with $\pmb{\psi}_i \in S^d : [0,1]^d\rightarrow\mathbb{R}^d$. The numerical approximation becomes, expressed in parametric coordinates, 
a vector-valued $d$-variate B-spline tensor product function (\ref{eq:bsplinef}).

\subsection{Numerical example}
We study now the magnitude of thermal expansion and its effect on the clearance height. To that end, it is enough to consider thermal growth of the rotors and the casing in a stationary case. We assume a known temperature distribution with the temperature growing linearly along the axial direction towards the high pressure section. In order to compare our results with~\cite{spille2015cfd}, we use the similar temperature distribution from 70$^\circ$C to 200$^\circ$C. The temperature at the shafts' end is also set to 70$^\circ$C. Both rotors are fixed at the low pressure end of the shaft, while the high pressure ends are allowed to move in the axial direction. In this setting both the deformation of the lateral surfaces and the axial elongation can be computed. The casing is fixed at the outer boundary and is allowed to expand inside. Its outer surface is set to 70$^\circ$C while the same linear temperature distribution is assumed at the inner surface. For the material parameters we choose those of steel: the thermal expansion coefficient $\alpha$ is 12e-6 K$^{-1}$, Young's modulus is 209 GPa and Poisson's ratio is 0.3. The computational meshes consist of 14592 cubic elements for the male rotor, of 18944 cubic elements for the female rotor and, after one iteration of uniform refinement, of 5184 quadratic elements for the casing. The simulation code was implemented using an open-source C++ library G+Smo, which provides most of the geometric routines needed in IGA~\cite{jlmmz2014}.

\begin{figure}[H]
	\centering
	\includegraphics[height = 4.51cm]{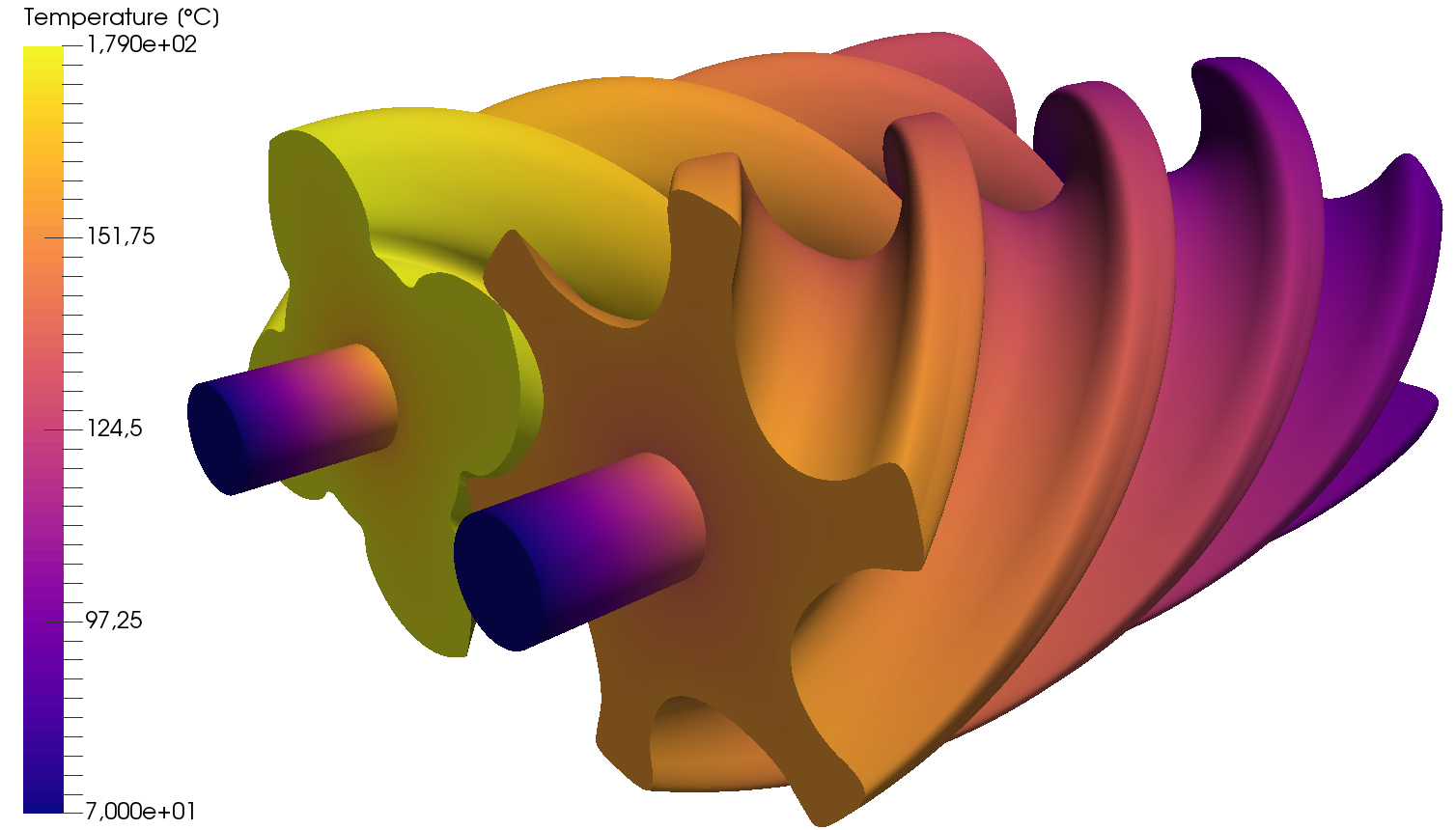}
	\includegraphics[height = 4.51cm]{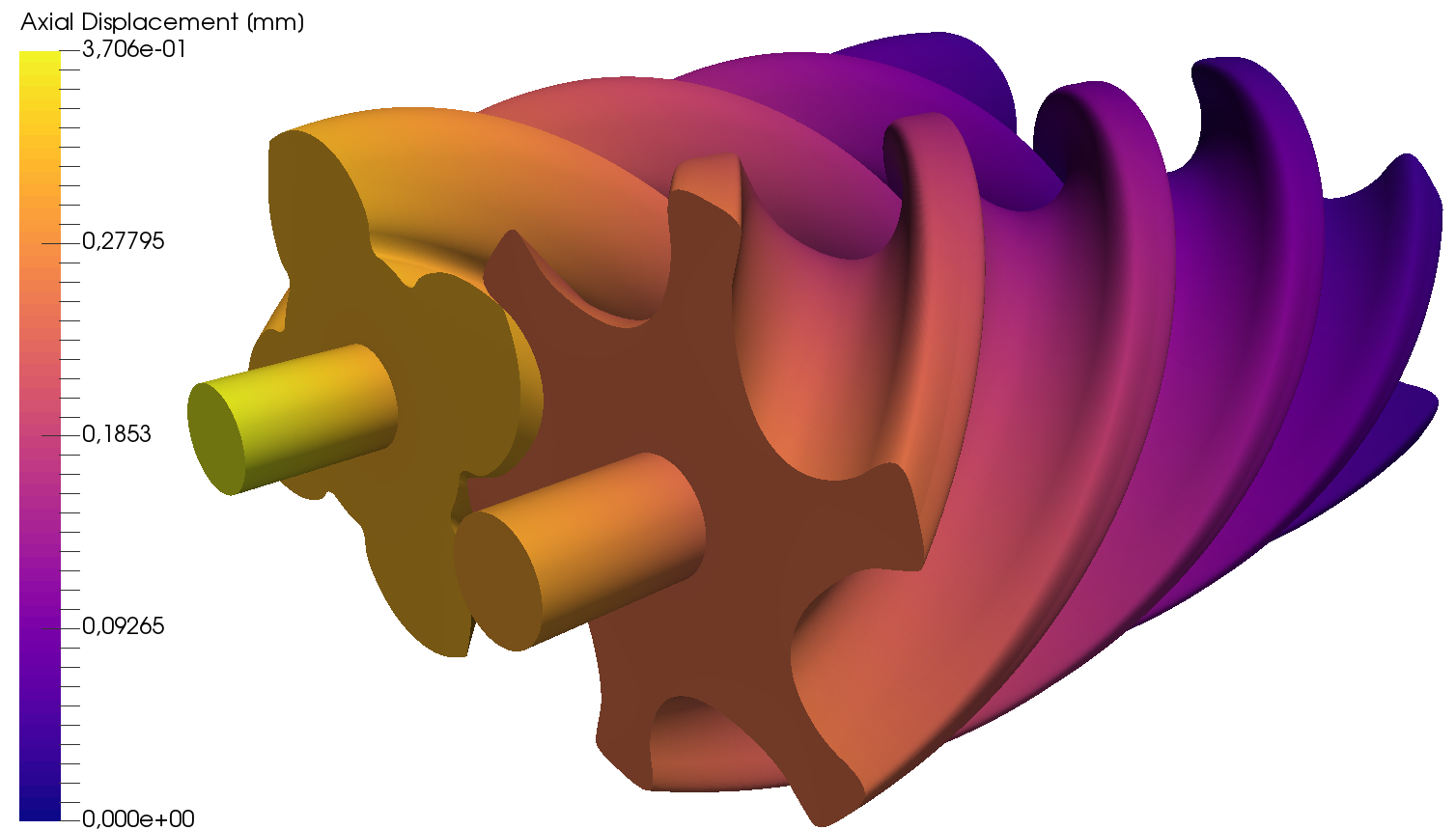}
	\caption{Temperature distribution for the compressor in operation (left). Axial elongation for the rotors (right).}\label{fig:temp}
\end{figure}

\begin{figure}[H]
	\centering
	\includegraphics[height = 5.2cm]{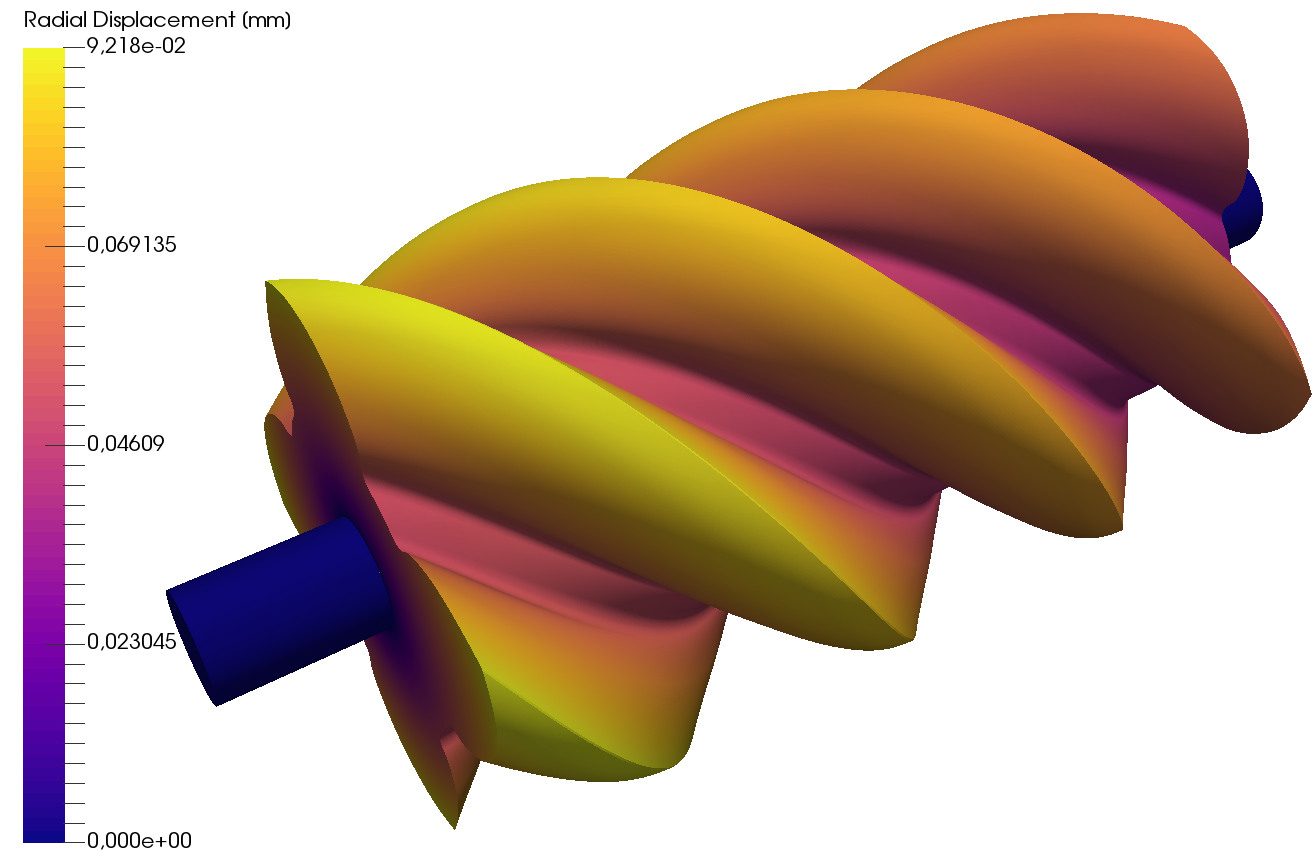}
	\includegraphics[height = 5.2cm]{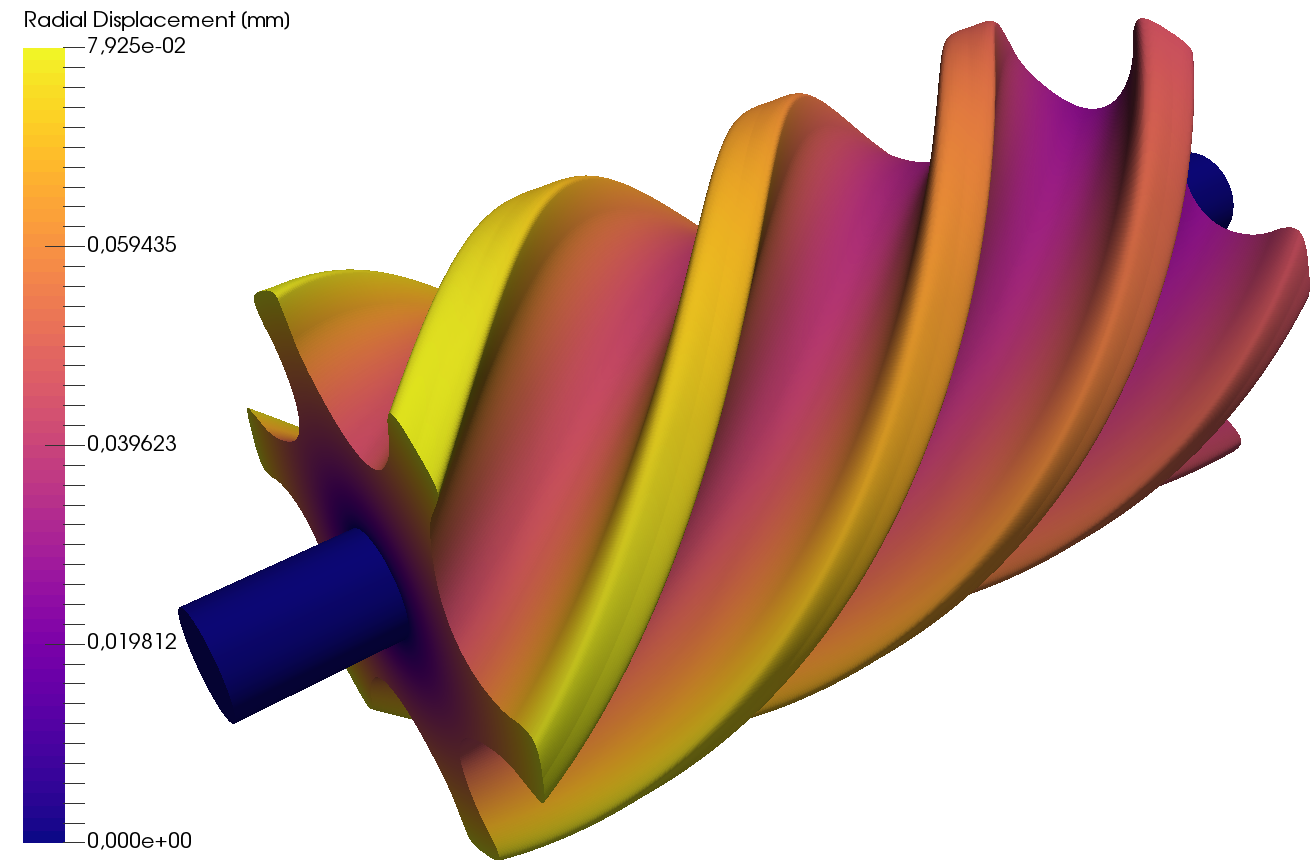}
	
	\caption{Radial expansion for the male rotor (left) and the female rotor (right).}\label{fig:disp}
\end{figure}

\begin{figure}[H]
	\centering
	\includegraphics[height = 5.11cm]{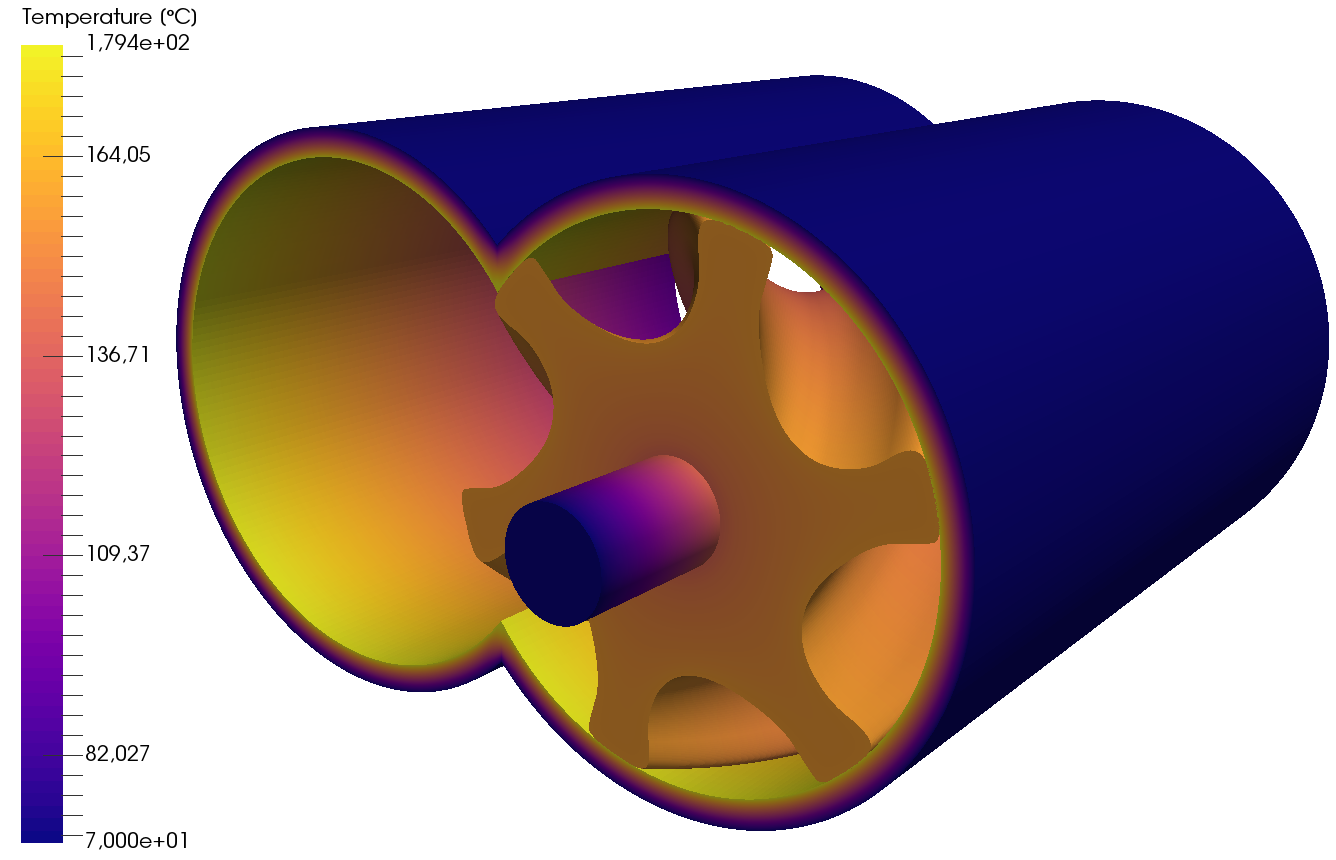}
	\includegraphics[height = 5.11cm]{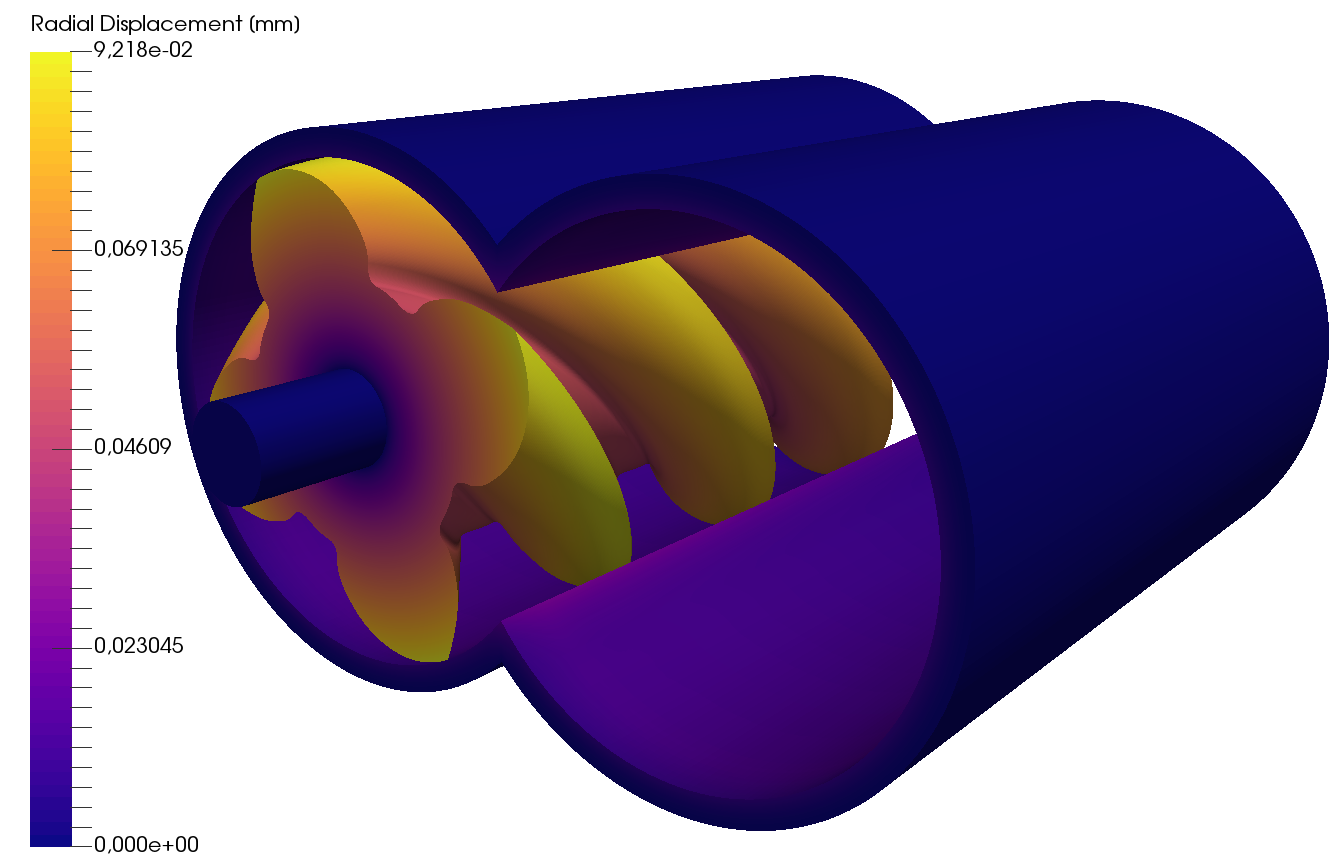}
	
	\caption{Simulation results with the casing. For the sake of visibility only one rotor is visualized. Temperature distribution (left) and radial expansion (right).}\label{fig:casing}
\end{figure}

\begin{figure}[H]
	\centering
	\includegraphics[height = 5.35cm]{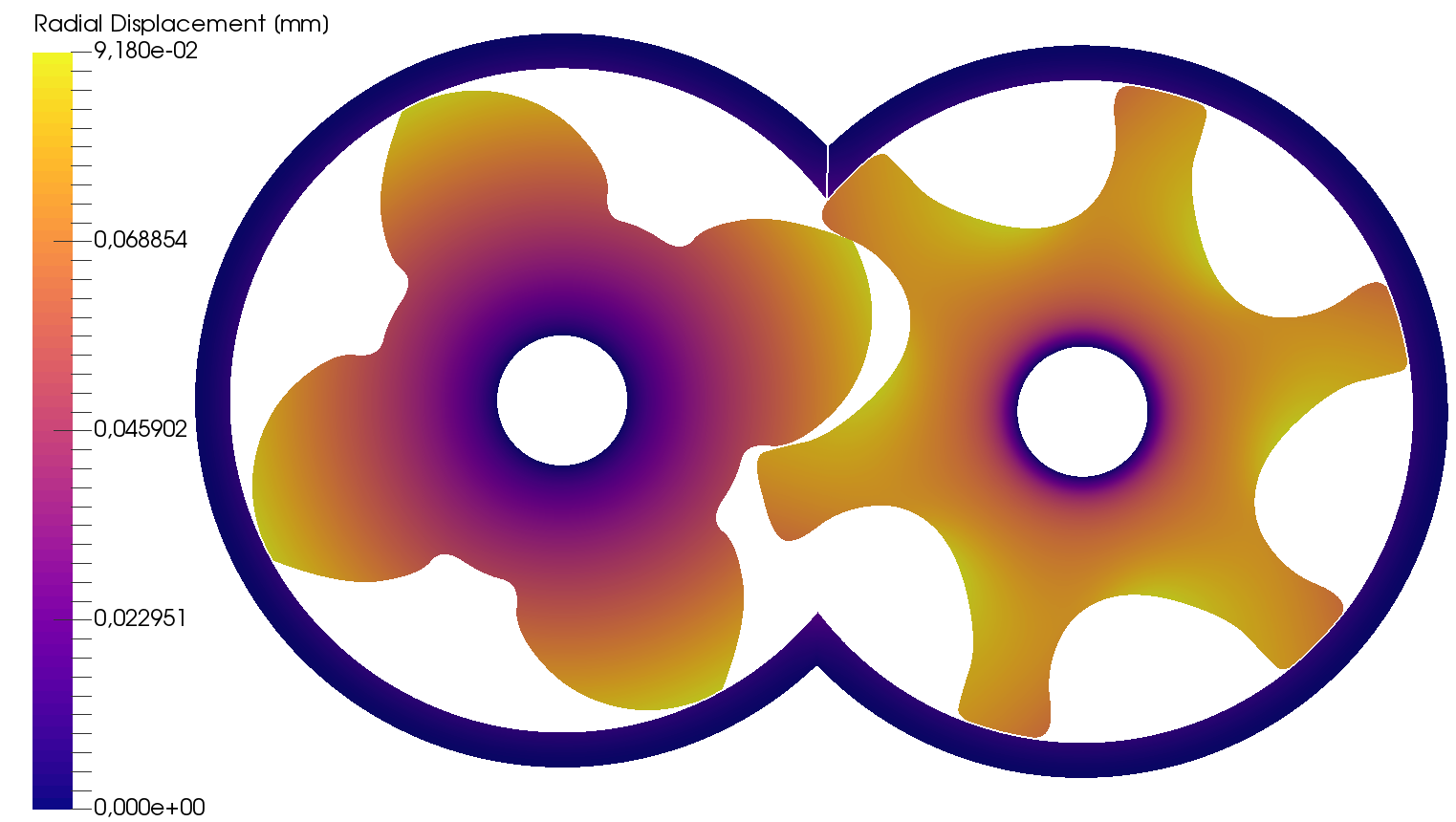}
	\includegraphics[height = 5.35cm]{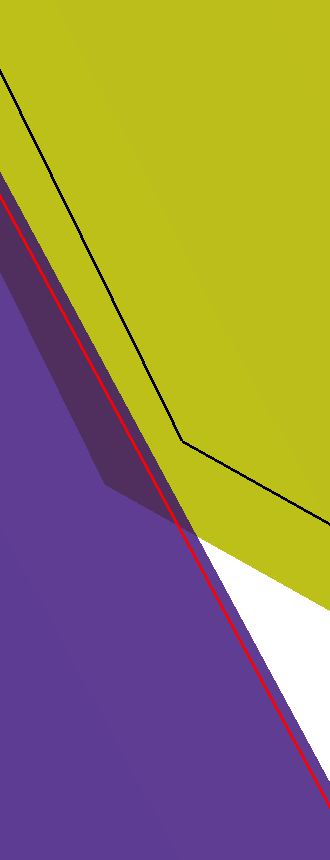}
	\includegraphics[height = 5.35cm]{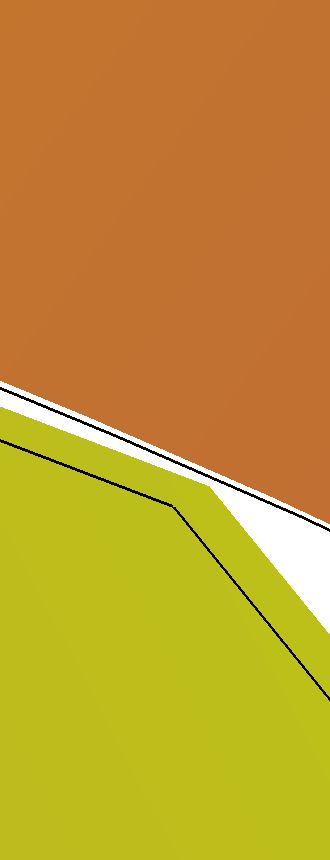}
	\includegraphics[height = 5.35cm]{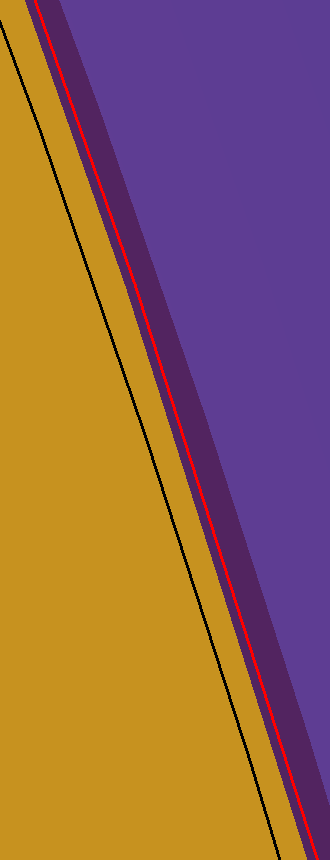}
	\caption{Clearance contraction. Axial expansion in 2D cross-section (left). Magnified views on the clearance space (right): between the male rotor and the casing (left); between the rotors (middle); between the female rotor and the casing (right).}\label{fig:dispZ}
\end{figure}

The results of the simulation are depicted in Fig.~\ref{fig:temp},~\ref{fig:disp} and~\ref{fig:casing}. The simulation predicts an axial elongation of 317 $\mu$m for the male rotor, 370 $\mu$m including the shaft. For the female rotor a predicted axial elongation is 251 $\mu$m, 300 $\mu$m including the shaft. Predicted radial expansion is 92 $\mu$m and 79 $\mu$m for the male and the female correspondingly. The casing expands for 12 $\mu$m towards rotors. The results on radial expansion are in a very good agreement with~\cite{spille2015cfd}. However, just like in the original paper they would mean a complete closing of the clearance between the casing and the rotors, which is originally 44 $\mu$m. Figure~\ref{fig:dispZ} shows a collection of magnified views on the clearance space in the 2D cross-section of the high pressure side of the rotors. Due to thermal expansion the clearance is completely closed between the casing and the rotors. Next to the sharp edge of the male the clearance space between the rotors is contracted by approximately 60\%.

\section{Conclusions and future work}
In this paper, we proposed to use Isogeometric Analysis to simulate thermal expansion of twin screw compressors. We demonstrated one possible approach to generating tensor-product B-spline parametrizations of the rotors and used them to solve a stationary thermal expansion problem. The results are consistent with the results from the previous works obtained with the help of the commercial software, but much less computational resources are necessary. However, further investigation is required to understand the excessive radial expansion which leads to a complete closing of the clearance space. In the future we plan to improve the model by using heat fluxes obtained from a coupled FSI simulation. This simulation will include the flow of the gas between the rotors and will take into account the deformation of the flow domain due to the thermal expansion. Furthermore, we plan to validate the simulation results on an experimental test rig with the help of our colleagues.  

\subsection*{Acknowledgement}
We are grateful to Andreas Brümmer and Matthias Utri for supplying us with the point cloud data of the rotors' geometry and for the fruitful discussions about the simulation setup. We are grateful to Matthias Möller and Jochen Hinz for valuable advice on the meshing strategy. This research is supported by the European Union under grant no. 678727 (Alexander Shamanskiy, Project MOTOR).

\bibliographystyle{iopart-num}
\bibliography{refs}

\end{document}